\documentclass{amsproc}
\usepackage{amsfonts}

\usepackage{graphicx}
\usepackage{amscd}
\usepackage{amsmath}

\theoremstyle{plain}

\newtheorem{corollary}{Corollary}

\newtheorem{example}{Example}

\newtheorem{proposition}{Proposition}
\newtheorem{remark}{Remark}

\newtheorem{theorem}{Theorem}
\numberwithin{equation}{section}

\begin{document}
\title[Automatic Continuity]{\vspace*{-0.4in}Automatic Continuity of Biseparating Maps}
\author{Jes\'us Araujo}
\address{Departamento de Matem\'{a}ticas, Estad\'{\i}stica y Computaci\'{o}n\\
Universidad de Cantabria\\
Facultad de Ciencias\\
Avda. de los Castros, s. n.\\
E-39071 Santander, Spain}
\email{araujoj@unican.es\qquad\ http://www.matesco.unican.es/\symbol{126}araujo/ }
\author{Krzysztof Jarosz}
\address{Department of Mathematics and Statistics\\
Southern Illinois University Edwardsville\\
IL 62026, USA}
\email{\ kjarosz@siue.edu\qquad http://www.siue.edu/\symbol{126}kjarosz/}
\thanks{Research of J. Araujo partially supported by the Spanish Direcci\'{o}n
General de Investigaci\'{o}n Cient\'{\i}fica y T\'{e}cnica (DGICYT,
PB98-1102)}

\begin{abstract}
We prove that a biseparating map between spaces of vector-valued continuous
functions is usually automatically continuous. However, we also discuss
special cases when it is not true.
\end{abstract}

\maketitle

\section{Introduction}

Assume $A,B$ are spaces of continuous functions on a set $X$ taking values
in a normed space $E$. A linear map $T:A\rightarrow B$ is called separating
if 
\begin{equation*}
\left\| f\left( \cdot \right) \right\| \left\| g\left( \cdot \right)
\right\| \equiv 0\qquad \Longrightarrow \qquad \left\| Tf\left( \cdot
\right) \right\| \left\| Tg\left( \cdot \right) \right\| \equiv 0\text{,}
\end{equation*}
and biseparating if $T^{-1}$ exists and is separating as well. The concept
has its source in the theory of topological lattices but is also an
important generalization of multiplicative maps on Banach algebras with
application to many other areas. It has been studied intensively by a number
of authors, see e.g. \cite{YA2000, AJ1, ABN1996, ABN, JA, Da, FH, JW}. The
most typical questions are whether in a given setting a separating linear
bijection must be biseparating, and if such a map must be automatically
continuous. Both questions have positive answers for separating linear
bijections between $C\left( X\right) $ spaces with $X$ compact \cite{KJ1},
negative in general \cite{YA}, and are open in some other cases.

In this note we show that any biseparating map between spaces of \emph{%
bounded} vector valued functions is automatically continuous if and only if
the underlying topological space $X$ has no isolated points. For maps
between spaces of unbounded functions the situation is more complicated but
perhaps also more interesting. If $X$ is metrizable, or connected, or
locally compact, or separable, or satisfies some other simple restrictions,
we again get the same condition - no isolated points. However, there are
sets $X$ with no isolated points that admit biseparating maps which are
discontinuous in a very strong sense. Our results are valid in both real and
complex cases.

\section{Preliminaries}

For a completely regular space $X$ and a normed space $E$ we denote by $%
C\left( X,E\right) $ and $C_{b}\left( X,E\right) $ the space of all
continuous functions from $X$ into $E$, and the space of all bounded
continuous functions from $X$ into $E$, respectively. If $E$ is equal to the
scalar field, we write $C\left( X\right) $ in place of $C\left( X,E\right) $%
. We equip $C_{b}\left( X,E\right) $ with the topology of uniform
convergence on $X$, and $C\left( X,E\right) $ with the topology of uniform
convergence on compact subsets of $X$. For normed spaces $E,F$ we denote by $%
L^{-1}\left( E,F\right) $ the set of all linear bijections from $E$ onto $F$.

For a function $f\in C\left( X\right) $ it may be often convenient to extend 
$f$ to a continuous function on $\beta X$, the maximal compactification of $%
X.$ In general, for $x\in \beta X\backslash X,$ the value of $f\left(
x\right) $ may be infinite. There are, however, completely regular spaces $X$
with a point $x_{0}\in \beta X\backslash X$ such that the value of $f\left(
x_{0}\right) $ is finite for all $f\in C\left( X\right) $; the set of all
points in $\beta X$ with this property is called the realcompactification of 
$X$ and is denoted by $rX.$ Since we have $C\left( X\right) =C\left(
rX\right) ,$ the natural domain for a continuous function on $X$ is $rX$,
not $X$. Hence we will often consider realcompactifications, or
alternatively, we will assume that the completely regular spaces under
consideration are realcompact. All compact sets are clearly realcompact,
also all subsets of Euclidean spaces are realcompact, or even more generally
all metrizable spaces of nonmeasurable cardinal are realcompact (\cite{GJ},
p.232).

We will often refer to the following result from \cite{AJ2}.

\begin{theorem}
\label{form}Suppose that $X$ and $Y$ are completely regular realcompact
spaces, $E,F$ are normed spaces, and that 

\begin{itemize}
\item  $T$ is a linear biseparating map from $C\left( X,E\right) $ onto $%
C\left( Y,F\right) $, 

\item  or $E,F$ are infinite-dimensional and $T$ is a linear biseparating
map from $C_{b}\left( X,E\right) $ onto $C_{b}\left( Y,F\right) $. 
\end{itemize}

Then there exists a homeomorphism $h:Y\rightarrow X$ and a map $%
J:Y\rightarrow L^{-1}(E,F)$ such that $(Tf)(y)=(J\left( y\right) )(f(h(y))$
for every $f$ in the domain of $T$, and $y\in Y$.
\end{theorem}

Notice the result can also be applied if $E,F$ are finite-dimensional and $T$
is defined on $C_{b}\left( X,E\right) $ since in such case we have $%
C_{b}\left( X,E\right) \simeq C_{b}\left( \beta X,E\right) =C\left( \beta
X,E\right) $. Based on the above Theorem, each time we are asked about
continuity of a biseparating map $T$, we can compose the map with $%
f\longmapsto f\circ h^{-1}$, and assume without loss of generality that $X=Y$
and that$\ $%
\begin{equation}
(Tf)(x)=(J\left( x\right) )(f(x)),\quad \text{ for every }f\text{ and }x\in
X.  \label{1}
\end{equation}
Hence, if the normed spaces $E,F$ are finite-dimensional any biseparating
map $T$ must be automatically continuous. On the other hand, if $X$ contains
an isolated point $x_{0}$ and $E$ is infinite-dimensional then one can
trivially define a discontinuous biseparating map on $C\left( X,E\right) $
or $C_{b}\left( X,E\right) $: take any discontinuous linear bijection $J_{0}$
from $E$ onto itself and put 
\begin{equation*}
\left( Tf\right) \left( x\right) =\left\{ 
\begin{array}{lll}
f\left( x\right)  & \text{if} & x\neq x_{0,} \\ 
J_{0}\left( f\left( x_{0}\right) \right)  & \text{if} & x=x_{0}.
\end{array}
\right. 
\end{equation*}
In the next section we show that, except for some very special sets $X$,
this is basically the only way one can get a discontinuous biseparating map.

Recall that a point $x_{0}$ in a topological space is called a $P$-point if
for each sequence $(U_{n})$ of neighborhoods of $x_{0}$, the intersection $%
\bigcap_{n=1}^{\infty}U_{n}$ contains again a neighborhood of $x_{0}$; if a
topological space is first countable, for example if it is metrizable, then
only isolated points are the $P$-points. A completely regular topological
space is called $P$-space if each point of the space is a $P$-point. Among
the metric spaces, $P$-spaces coincide with the discrete spaces. We will
need the following properties of the $P$-spaces (\cite{GJ}, p. 63).

\begin{proposition}
\label{Prop GJ}For a completely regular space $X$ the following conditions
are equivalent

\begin{itemize}
\item  $X$ is a $P$-space,

\item  any real-valued continuous function on $X$ is locally constant,

\item  any $G_{\delta}$ subset of $X$ is open.
\end{itemize}
\end{proposition}

\section{Results}

\begin{theorem}
\label{th1}Assume $X,Y$ are completely regular realcompact spaces, $E,F$ are
normed spaces, and $T:C_{b}(X,E)\rightarrow C_{b}(Y,F)$ is a biseparating
map. If $Y$ has no isolated points then $T$ is continuous.
\end{theorem}

\begin{proof}
By Theorem \ref{form} we may assume without loss of generality that $Y=X$
and that $T$ is of the form (\ref{1}). Suppose that $T$ is not continuous
and let $(f_{n})_{n=1}^{\infty }$ be a sequence of elements in $C_{b}(X,E)$
such that 
\begin{equation*}
\left\| f_{n}\right\| <1/n^{2}\text{ and }\left\| Tf_{n}\right\|
>n+\sum_{j=1}^{n-1}\left\| Tf_{j}\right\| \text{, \quad for }n\in \mathbb{N}.
\end{equation*}
Let $x_{1}\in X$ be such that $\left\| (Tf_{1})(x_{1})\right\| \geq 1$.
Since $x_{1}$ is not isolated, then there exists an $x_{2}\neq x_{1}$ with $%
\left\| (Tf_{2})(x_{2})\right\| \geq 1+\left\| Tf_{1}\right\| $. Continuing
the process, we can inductively construct a sequence $\left( x_{n}\right)
_{n=1}^{\infty }$ of different points in $X$ such that 
\begin{equation*}
\left\| (Tf_{n})(x_{n})\right\| \geq n+\sum_{j=1}^{n-1}\left\|
Tf_{j}\right\| \text{, \quad for }n\in \mathbb{N}.
\end{equation*}
Let $(U_{n})_{n=1}^{\infty }$ be a sequence of open subsets of $X$ and $%
(g_{n})$ a sequence of continuous bounded real-valued functions on $X$ such
that for any $n\in \mathbb{N}$%
\begin{equation*}
x_{n}\in U_{n}\text{ and }x_{k}\notin U_{n}\text{, \quad for }k<n,
\end{equation*}
and 
\begin{equation*}
g_{n}\left( x_{n}\right) =1=\left\| g_{n}\right\| \quad \text{and}\quad 
\mathrm{supp}\;g_{n}\subset U_{n}.
\end{equation*}

Since $\left\| g_{n}f_{n}\right\| \leq \left\| f_{n}\right\| <1/n^{2}$ we
can define a function $f$ in $C_{b}(X,E)$ by 
\begin{equation*}
f:=\sum_{j\in \mathbb{N}}g_{j}f_{j}.
\end{equation*}
>From (\ref{1}) and since $g_{n}\left( x_{n}\right) =1$ and $g_{j}\left(
x_{n}\right) =0$ for $j>n$ we get 
\begin{align*}
\left\| \left( Tf\right) \left( x_{n}\right) \right\| & =\left\| (J\left(
x_{n}\right) )(f(x_{n}))\right\| =\left\| (J\left( x_{n}\right) )\left(
\sum_{j\in \mathbb{N}}g_{j}\left( x_{n}\right) f_{j}\left( x_{n}\right)
\right) \right\|  \\
& =\left\| (J\left( x_{n}\right) )\left( \sum_{j=1}^{n-1}g_{j}\left(
x_{n}\right) f_{j}\left( x_{n}\right) +f_{n}\left( x_{n}\right) \right)
\right\|  \\
& =\left\| \sum_{j=1}^{n-1}g_{j}\left( x_{n}\right) (J\left( x_{n}\right)
\left( f_{j}\left( x_{n}\right) \right) +(J\left( x_{n}\right) )\left(
f_{n}\left( x_{n}\right) \right) \right\|  \\
& =\left\| \sum_{j=1}^{n-1}g_{j}\left( x_{n}\right) (T\left( f_{j}\right)
\left( x_{n}\right) +(T\left( f_{n}\right) )\left( x_{n}\right) \right\|  \\
& \geq \left\| (T\left( f_{n}\right) )\left( x_{n}\right) \right\|
-\sum_{j=1}^{n-1}\left\| (T\left( f_{j}\right) \left( x_{n}\right) \right\| 
\\
& \geq \left( n+\sum_{j=1}^{n-1}\left\| Tf_{j}\right\| \right)
-\sum_{j=1}^{n-1}\left\| Tf_{j}\right\| \geq n.
\end{align*}
Since our space contains only bounded functions, the contradiction shows
continuity of $T$.
\end{proof}

\begin{theorem}
\label{th2}Assume $X,Y$ are completely regular realcompact spaces and $E,F$
are both infinite-dimensional normed spaces. Then

\begin{itemize}
\item[\textbf{(i) \ }]  if the interior of the set of $P$-points of $Y$ is
empty then every biseparating map $T:C(X,E)\rightarrow C(Y,F)$ is continuous,

\item[\textbf{(ii) }]  if the interior of the set of $P$-points of $Y$ is
not empty then there is a discontinuous biseparating map $%
T:C(Y,F)\rightarrow C(Y,F),$

\item[\textbf{(iii)}]  if $\dim E=\dim F$ then all points of $Y$ are $P$%
-points if and only if for any linear bijection $J$ from $E$ onto $F$
(continuous or not) the map $T\left( f\right) \overset{df}{=}J\circ f$ is
well defined and biseparating from $C(Y,E)$ onto $C(Y,F)$.
\end{itemize}
\end{theorem}

\begin{corollary}
For completely regular realcompact spaces $X,Y$ and normed spaces $E,F$
there is a biseparating map between $C(X,E)$ and $C(Y,F)$ if and only if $X$
and $Y$ are homeomorphic and

\begin{itemize}
\item  $E$ and $F$ are isomorphic, or

\item  $\dim E=\dim F$ and $Y$ is a $P$-space.
\end{itemize}
\end{corollary}

\begin{remark}
There is another useful way to think about these results. Notice that if $J$
is a linear bijection from a normed space $E$ onto itself, then we can
phrase the usual definition of continuity of $J$ as follows:

$J$\emph{\ is continuous iff transformation }$f\longmapsto J\left( f\right) $%
\emph{\ maps }$C\left( \mathbb{N}^{\ast},E\right) $\emph{\ into itself,}

where by $\mathbb{N}^{\ast }$ we denote the one-point compactification of
the set of natural numbers. Now, based on our results we can replace the set 
$\mathbb{N}^{\ast }$ in the definition above with a set $X$ if and only if $%
X $ is not a $P$-space.
\end{remark}

\begin{proof}[Proof of the Theorem]
\textbf{(i) }Assume that the interior of the set of $P$-points of $Y$ is
empty and $T:C(X,E)\rightarrow C(Y,F)$ is a biseparating map. As before by
Theorem \ref{form} we may assume without loss of generality that $Y=X$ and
that $T$ is of the form (\ref{1}).

First we are going to prove that $J(y):E\rightarrow F$ is continuous for
every $y\in X$. Assume that there is a $y_{0}\in X$ such that $J\left(
y_{0}\right) $ is discontinuous. If $y_{0}$ is not a $P$-point we put $x_{0}%
\overset{df}{=}y_{0}$. If $y_{0}$ is a $P$-point we put $V_{n}=\left\{ x\in
X:\left\| J\left( x\right) \right\| >n\right\} $. Since the map $\left\|
\cdot \right\| \circ J:X\rightarrow \lbrack 0,+\infty ]$ sending each $y\in
X $ into $\left\| J(y)\right\| $ is continuous, all of the sets $V_{n}$ are
neighborhoods of $y_{0}$, and so is $\bigcap_{n=1}^{\infty }V_{n}$. Moreover
for any $x\in \bigcap_{n=1}^{\infty }V_{n}$ the map $J\left( x\right) $ is
discontinuous. Since the interior of the set of $P$-points is empty, the set 
$\bigcap_{n=1}^{\infty }V_{n}$ contains a point which is not a $P$-point. We
denote by $x_{0}$ one of such points.

Let $(e_{n})$ be a sequence in $E$ such that 
\begin{equation*}
\left\| e_{n}\right\| <1/n^{2}\text{ and }\left\| J\left( x_{0}\right)
\left( e_{n}\right) \right\| >n\text{, \quad for }n\in\mathbb{N}.
\end{equation*}
Let $\mathbf{e}_{n}$ be the constant function equal to $e_{n}.$ Since $x_{0}$
is not a $P$-point there is a sequence of open neighborhoods $U_{n}$ of $%
x_{0}$ with 
\begin{equation*}
U_{n}\subset\{x\in X:\left\| (T\mathbf{e}_{n})(x)\right\| =\left\| J\left(
x\right) \left( e_{n}\right) \right\| >n\}
\end{equation*}
and such that the intersection of all $U_{n}$ is not a neighborhood of $%
x_{0} $; we may also assume that $\overline{U_{n+1}}\subset U_{n}$, for $n\in%
\mathbb{N}$. It follows that $x_{0}$ belongs to the closure of the following
set: 
\begin{equation*}
A:= \bigcup_{n=0}^{\infty}U_{n+1}-U_{n+2}.
\end{equation*}

For each $n\in\mathbb{N}$, let $g_{n}\in C_{b}(X)$ be such that 
\begin{equation*}
\left\| g_{n}\right\| =1=g_{n}\left( x\right) \text{, for }x\in U_{n+1}\text{%
, and }\mathrm{supp}\;g_{n}\subset U_{n}.
\end{equation*}

Put 
\begin{equation*}
g:=\sum_{n=1}^{\infty}(g_{n}-g_{n+2})e_{n+1};
\end{equation*}
since $\left\| g_{n}-g_{n+2}\right\| \leq2$ and $\left\| e_{n}\right\|
<1/n^{2}$ the series is convergent and $g\in C(X,E)$. For any $k\in\mathbb{N}
$ and any $x_{k+1}\in U_{k+1}-U_{k+2}$ we have 
\begin{equation*}
(g_{n}-g_{n+2})\left( x_{k+1}\right) =\left\{ 
\begin{array}{lll}
1-1=0 & \text{if} & n<k \\ 
1-0=1 & \text{if} & n=k \\ 
0-0=0 & \text{if} & n>k
\end{array}
\right. .
\end{equation*}
Hence 
\begin{align*}
\left\| (Tg) (x_{k+1})\right\| & =\left\| \left( J\left( x_{k+1}\right)
\right) \left( \sum_{n=1}^{\infty}(g_{n}-g_{n+2})\left( x_{k+1}\right)
e_{n+1}\right) \right\| \\
& =\left\| \left( J\left( x_{k+1}\right) \right) \left( e_{k+1}\right)
\right\| =\left\| (T\mathbf{e}_{k+1})(x_{k+1})\right\| >k+1.
\end{align*}

Since $x_{0}$ belongs to the closure of $\bigcup_{n=0}^{%
\infty}U_{n+1}-U_{n+2}$ we get $\left\| (Tg) (x_0)\right\| =\infty$, and
this contradiction proves the continuity of $J(y)$ for every $y \in X$.

Suppose now that $T$ is not continuous. Then there is a compact subset $K$
of $Y$ and a sequence $(f_{n})$ in $C(X,E)$ with $p_{K}(f_{n})<1/n^{2}$ and $%
p_{K}(Tf_{n})>n$.

Let $y_{n}\in K$ be such that $\left\| (Tf_{n})(y_{n})\right\| >n,$ for each 
$n\in \mathbb{N}$. Notice that the set $\left\{ y_{n}\in K:n\in \mathbb{N}%
\right\} $ is infinite since otherwise there would be a point $y_{n_{0}}\in
K $ with $\left\| (Tf_{n})(y_{n_{0}})\right\| >n$ for infinitely many $n\in 
\mathbb{N}$, contrary to the continuity of $J\left( y_{n_{0}}\right) $.
Without loss of generality we may assume that $y_{n}\neq y_{m}$ when $n\neq
m $, and take a sequence $(U_{n})_{n=1}^{\infty }$ of pairwise disjoint open
subsets of $X$ such that $y_{n}\in U_{n}\subset \left\{ x\in X:\left\|
f_{n}\left( x\right) \right\| <1/n^{2}\right\} ,$ for $n\in \mathbb{N}$. Let 
$(g_{n})_{n=1}^{\infty }$ be a sequence in $C_{b}(X)$ such that $\mathrm{supp%
}\;g_{n}\subset U_{n}$ and $g_{n}(y_{n})=1,$ for $n\in \mathbb{N}$. It is
easy to see that 
\begin{equation*}
f\overset{df}{=}\sum_{k=1}^{\infty }g_{k}f_{k}\in C(X,E)
\end{equation*}
while, by \ref{1}, $\left\| (Tf)(y_{n})\right\| =\left\|
(Tf_{n})(y_{n})\right\| >n$ so $Tf$ is not bounded on $K$. Since this is not
possible we conclude that $T$ must be continuous.

\textbf{(ii) \ }Assume now that the interior $\mathcal{I}$ of the set of $P$%
-points of $Y$ is not empty and let $h_{0}$ be a continuous real-valued
function on $Y$ such that $h_{0}\left( x\right) =1$ for $x\in X\backslash 
\mathcal{I}$, and $h_{0}\left( z_{0}\right) =0$ at some point $z_{0}\in 
\mathcal{I}$. Put $\mathcal{A}=h_{0}^{-1}\left( 0\right) $. By Proposition
\nolinebreak \ref{Prop GJ} $\mathcal{A}$ is an open subset of $\mathcal{I}$,
and consequently of $X$; clearly as a zero set $\mathcal{A}$ is also closed.
Let $J$ be a discontinuous linear bijection from $F$ onto itself and define $%
T:C(Y,F)\rightarrow C(Y,F)$ by 
\begin{equation*}
\left( Tf\right) \left( x\right) =\left\{ 
\begin{array}{lll}
f\left( x\right) & \text{if} & x\notin \mathcal{A} \\ 
J\left( f\left( x\right) \right) & \text{if} & x\in \mathcal{A}
\end{array}
.\right.
\end{equation*}
The map $T$ is obviously biseparating but we need to check that it is well
defined, that is, that $Tf$ is continuous for any continuous $f$. It is of
course true on the clopen subset $X\backslash \mathcal{A}$. Let $f\in C(Y,F)$
and let $z\in \mathcal{A}$. By Proposition \ref{Prop GJ} the function 
\begin{equation*}
\mathcal{A}\ni x\longmapsto \left\| f\left( x\right) -f\left( z\right)
\right\|
\end{equation*}
is constant on a neighborhood of $z$, hence $f$, and so $Tf,$ are both
locally constant and consequently continuous.

\textbf{(iii) }As above we assume without loss of generality that $X=Y$ and
that $h$ is the identity map. Suppose that $\dim E=\dim F$ and that there is
a point $x_0 \in X$ which is not a $P$-point. We well see that there exists
a linear bijection $J$ from $E$ onto $F$ such that the map $T\left( f\right) 
\overset{df}{=}J\circ f$ is not well defined.

Since $J$ is not continuous, then there exists a sequence $(e_n)$ of
norm-one elements of $E$ such that $\left\| J e_n \right\|>n^3$, for $n \in 
\mathbb{N}$. On the other hand, since $x_0$ is not a $P$-point, there exists
a sequence $(U_n)_{n=1}^{\infty}$ of neighborhoods of $x_0$ such that the
intersection of all $U_n$ is not a neighborhood of $x_0$. We may assume that 
$\overline{U_{n+1}}\subset U_{n}$, for $n \in \mathbb{N}$. Consequently $x_0$
belongs to the closure of at least one of the following sets: 
\begin{equation*}
A_{1}:=\bigcup_{n=1}^{\infty}U_{2n}-U_{2n+2}\quad\text{or}\quad
A_{2}:=\bigcup_{n=1}^{\infty}U_{2n-1}-U_{2n+1}.
\end{equation*}
Without loss of generality we may assume that $x_{0}\in\overline{A_{2}}$.

For each $n\in\mathbb{N}$, let $g_{n}\in C_{b}(X)$, $0 \le g_n \le 1$, be
such that 
\begin{equation*}
g_{n} (U_{2n-1}- U_{2n +1}) =1 \text{ and }\mathrm{supp}\;g_{n}\subset
U_{2n-2}- U_{2n+2}.
\end{equation*}

It is clear that the map 
\begin{equation*}
f\overset{df}{=}\sum_{k=1}^{\infty} \frac{1}{k^2}g_{k}f_{k}
\end{equation*}
belongs to $C(X,E)$. Also, according to the description of $T$ and to the
fact that $x_0 \in \overline{A_{2}}$, we deduce that $\left\| (Tf) (x_0)
\right\| = \infty$, which is absurd. Consequently, $X$ is a $P$-space.

To prove the other implication we just repeat the above argument in \textbf{%
(ii)} with $\mathcal{A}=X$.
\end{proof}

As for the proof of the Corollary, notice that if $Y$ is not a $P$-space,
then it contains a point $y_0$ which is not a $P$-point, and consequently,
as in the proof of Theorem~\ref{th2} \textbf{(i)}, $J(y_0)$ is continuous
(and bijective). On the other hand, by Theorem~\ref{form}, $(T^{-1} g) (x) =
(Kx) (g(h^{-1}(x)))$, for some $K \in L^{-1} (F,E)$. In particular $K (h(y))$
is the inverse of $J(y)$ (see the proof of Theorem 2.4 in \cite{AJ2}), which
implies in particular that $h(y)$ is not a $P$-point and, consequently, $E$
and $F$ are isomorphic.

\bigskip

The $P$-sets without isolated points are rather unusual, they can not be
locally compact or metrizable. An example of such a set can be found in \cite
{GJ}, below we provide another one.

\begin{example}
Let $\Omega $ be an uncountable set and put $X=\left\{ 0,1\right\} ^{\Omega
} $. We say that a function $f$ on $X$ depends on countably many coordinates
if there is a countable subset $\Omega _{0}$ of $\Omega $ such that 
\begin{equation*}
f\left( x\right) =f\left( x^{\prime }\right) \text{ whenever }x\left( \omega
\right) =x^{\prime }\left( \omega \right) \text{ for all }\omega \in \Omega
_{0}.
\end{equation*}
Notice that $f$ depends on countable many coordinates if and only if $%
f=f_{0}\circ p_{\Omega _{0}}$\ where $\Omega _{0}$ is a countable subset of $%
\Omega $, $p_{\Omega _{0}}$ is a natural projection from $\left\{
0,1\right\} ^{\Omega }$ onto $\left\{ 0,1\right\} ^{\Omega _{0}}$ and $f_{0}$
is a function on $\left\{ 0,1\right\} ^{\Omega _{0}}$. We denote by $%
C_{\aleph _{0}}\left( X\right) $ and $C_{\aleph _{0}}^{b}\left( X\right) $
the set of all scalar valued functions that depend on countably many
variables, and bounded functions with that property, respectively. We
introduce on $X$ the coarsest topology such that all functions in $C_{\aleph
_{0}}\left( X\right) $ are continuous. That topology $\tau $ of $X$ has a
base that can be described as follows 
\begin{equation*}
\left\{ p_{\Omega _{0}}^{-1}\left( A\right) :A\subset \left\{ 0,1\right\}
^{\Omega _{0}},card\left( \Omega _{0}\right) \leq \aleph _{0}\right\} .
\end{equation*}
Since a uniform limit of a sequence of functions depending on countably many
coordinates each, has also such a property, $C_{\aleph _{0}}^{b}\left(
X\right) $ is a self-adjoint uniform algebra so $X$ can be seen as a subset
of the maximal ideal space of that algebra. Since any subset of a compact
Hausdorff space is completely regular, so is $X$.

Suppose that $f$ is a continuous real-valued function on $X$ and fix $%
x_{0}\in X$. Without loss of generality, we assume that $x_{0}(\omega )=0$,
for $\omega \in \Gamma \subset \Omega $, and $x_{0}(\omega )=1$, for $\omega
\in \Omega \setminus \Gamma $. Since, for any $n\in N,$ the set $%
U_{n}=\left\{ x\in X:\left| f\left( x\right) -f\left( x_{0}\right) \right| <%
\frac{1}{n}\right\} $ is open, it contains an open set of the form 
\begin{equation*}
V_{n}=\left\{ x\in X:x\left( \Omega _{1,n}\right) \subset \left\{ 0\right\}
,x\left( \Omega _{2,n}\right) \subset \left\{ 1\right\} \right\} ,
\end{equation*}
for some countable sets $\Omega _{1,n}\subset \Gamma ,\Omega _{2,n}\subset
\Omega \setminus \Gamma $. Put $\Omega _{1}=\bigcup_{n=1}^{\infty }\Omega
_{1,n},\Omega _{2}=\bigcup_{n=1}^{\infty }\Omega _{2,n}$, clearly $\Omega
_{1},\Omega _{2}$ are countable and disjoint. Now 
\begin{equation*}
U_{0}=\left\{ x\in X:x\left( \Omega _{1}\right) \subset \left\{ 0\right\}
,x\left( \Omega _{2}\right) \subset \left\{ 1\right\} \right\} ,
\end{equation*}
is a neighborhood of $x_{0}$ and $f$ is constant on $U_{0}$.
\end{example}

\end{document}